\documentclass[11pt,a4paper]{article}

\usepackage{amsmath,amsthm,amssymb,amsfonts}
\usepackage{mathtools}
\usepackage[utf8]{inputenc}
\usepackage[T1]{fontenc}
\usepackage{hyperref}
\usepackage{enumitem}
\usepackage[margin=1in]{geometry}
\theoremstyle{plain}
\newtheorem{theorem}{Theorem}[section]
\newtheorem{proposition}[theorem]{Proposition}

\newtheorem{corollary}[theorem]{Corollary}
\theoremstyle{definition}
\newtheorem{definition}[theorem]{Definition}
\newtheorem{example}[theorem]{Example}
\theoremstyle{remark}
\newtheorem{remark}[theorem]{Remark}

\newcommand{\E}{\mathbb{E}}
\newcommand{\R}{\mathbb{R}}
\newcommand{\Id}{\mathrm{Id}}
\newcommand{\Var}{\mathrm{Var}}
\newcommand{\dom}{\mathrm{dom}}
\newcommand{\op}{\mathrm{op}}
\newcommand{\Mall}{\mathrm{Mall}}

\title{Stochastic Calculus as Operator Factorization:\\
An Operator--Covariant Derivative and Unified Clark--Ocone Representation}
\author{Ramiro Fontes\\
Quijotic Research\\
\texttt{ramirofontes@gmail.com}}
\date{February 2026}

\begin{document}
\maketitle

\begin{abstract}
We present a unified operator-theoretic framework for stochastic calculus based on the factorization
\[
(\Id - \E)\,F \;=\; \delta_X\,\Pi_X\,D_X\,F,
\]
valid for $\mathcal{F}_T^X$-measurable $F \in L^2(\Omega)$ when the driving process~$X$ has the representation property. This decomposes fluctuations from the mean into three operations determined by the chosen stochastic integral and filtration: differentiation, predictable projection, and integration.

For a square-integrable process~$X$ with stochastic integral $\delta_X \colon \mathcal{H}_X \to L^2(\Omega)$, we define the \emph{operator--covariant derivative} $D_X := \delta_X^*$ as the Hilbert space adjoint of~$\delta_X$. Combined with predictable projection~$\Pi_X$, this yields a unified Clark--Ocone representation
\[
F \;=\; \E[F] + \int_0^T (\Pi_X D_X F)_t\, dX_t
\]
for $\mathcal{F}_T^X$-measurable $F \in L^2(\Omega)$. The operator $D_X F$ is defined as an adjoint for all $F \in L^2(\Omega)$, without differentiability assumptions; the representation holds when $X$ has the predictable representation property, and reduces to the Galtchouk--Kunita--Watanabe projection (optimal $L^2$-approximation by stochastic integrals) when it does not.

The framework requires no reproducing kernel Hilbert space or Cameron--Martin structure, and applies to non-Gaussian processes. We work out concrete examples including Brownian motion (recovering the classical Clark--Ocone formula), general continuous martingales, and compensated Poisson processes, where we compute the operator--covariant derivative explicitly and verify the factorization.


\end{abstract}

%=============================================================================
\section{Introduction}
\label{sec:intro}
%=============================================================================

Two central themes run through modern stochastic analysis: stochastic integral representations (martingale representation, Clark--Ocone) and change-of-variables formulas (It\^o and its functional extensions). Although often developed as separate frameworks, both express a common algebraic structure that can be written as a single operator factorization:
\begin{equation}\label{eq:factorization}
(\Id - \E)\,F \;=\; \delta_X\,\Pi_X\,D_X\,F,
\end{equation}
valid for $F$ in the natural domain (square-integrable, $\mathcal{F}_T^X$-measurable) when the driving process has the representation property. Here $\Id$ is the identity, $\E$ is expectation, $D_X$ is a stochastic derivative, $\Pi_X$ is orthogonal projection onto predictable integrands, and $\delta_X$ is the stochastic integral (divergence) operator.

This paper develops the viewpoint that a significant part of stochastic calculus can be organized as operator factorization. The key construction is to define the stochastic derivative as the adjoint of the stochastic integral in the energy space of the driving process. For a square-integrable process~$X$ with bounded stochastic integral $\delta_X \colon \mathcal{H}_X \to L^2(\Omega)$, we set
\[
D_X := \delta_X^*
\]
via the Riesz representation theorem. Combined with the predictable projection~$\Pi_X$, this yields the unified Clark--Ocone representation
\begin{equation}\label{eq:CO}
F \;=\; \E[F] + \delta_X(\Pi_X D_X F)
\;=\; \E[F] + \int_0^T (\Pi_X D_X F)_t\, dX_t
\end{equation}
for $\mathcal{F}_T^X$-measurable $F \in L^2(\Omega)$, requiring only square integrability and representability of~$F$ rather than Sobolev-type differentiability.

%--- Novelty discussion ---
\subsection{Relationship to existing work and novelty}
\label{ssec:novelty}

The individual ingredients of our construction are well established. The idea that a stochastic derivative can be defined as the adjoint of the divergence operator is a foundational principle of Malliavin calculus~\cite{Nualart2006}. The Clark--Ocone formula is classical~\cite{Nualart2006,Ocone1984}. The martingale representation property has been studied extensively~\cite{Revuz1999}.

\paragraph{The circularity question.}
We address directly a natural objection: Theorem~\ref{thm:CO} assumes the representation property (PRP) and concludes a representation formula, so one might ask whether the result is tautological. The content of the theorem is not the \emph{existence} of a representation integrand---that is assumed---but rather its \emph{identification} as $\Pi_X D_X F$, i.e., the predictable projection of the adjoint of $F$ against the stochastic integral. This identification is what makes the factorization~\eqref{eq:factorization} an operator identity rather than merely an existence statement. In the Brownian case, this recovers the known formula $(D_B F)_t = \E[D_t^{\Mall} F \mid \mathcal{F}_t]$ without passing through the Malliavin--Sobolev space $\mathbb{D}^{1,2}$; in the Poisson case, it recovers the predictable projection of the add-a-point operator (see Proposition~\ref{prop:poisson_explicit}). These are not new results, but the operator-algebraic viewpoint organizes them under a single identity.

\paragraph{What is new.}
The contribution is one of organization and perspective, not of new hard analysis. Specifically:

\begin{enumerate}[label=(\arabic*),leftmargin=2em]
\item A single operator identity~\eqref{eq:factorization} unifies Clark--Ocone representations across different classes of driving processes---including non-Gaussian ones---without requiring the reproducing kernel Hilbert space (RKHS) or Cameron--Martin structures on which classical Malliavin calculus depends. The hypotheses (centered integrals, isometry on predictable processes, PRP) are transparent and can be checked process by process.

\item The operator $D_X = \delta_X^*$ is defined on all of $L^2(\Omega)$ by Riesz representation, not merely on a dense Sobolev-type subspace such as the Malliavin--Sobolev space $\mathbb{D}^{1,2}$~\cite{Nualart2006}. This full-domain property has concrete content: when PRP fails, $\delta_X(\Pi_X D_X F)$ is the Galtchouk--Kunita--Watanabe projection of $F - \E[F]$ onto the stable subspace generated by~$X$ (Proposition~\ref{prop:GKW}), i.e., the variance-minimizing approximation by stochastic integrals. The factorization~\eqref{eq:factorization} is thus an exact identity under PRP and a best-approximation formula without it.

\item We verify the framework against a substantive non-Gaussian example: the compensated Poisson process. In Section~\ref{ssec:poisson}, we compute $D_{\tilde{N}}$ for Poisson functionals and confirm that it recovers the well-known finite-difference (add-a-point) operator from the discrete Malliavin calculus of~\cite{Privault2009,DiNunno2009}. The novelty is not the operator itself but its characterization as an instance of the general adjoint construction.
\end{enumerate}

\subsection{Organization}
Section~\ref{sec:setup} establishes notation. Section~\ref{sec:OCD} introduces the operator--covariant derivative and proves the unified Clark--Ocone representation. Section~\ref{sec:examples} treats concrete examples: Brownian motion, general martingales, and compensated Poisson processes. Section~\ref{sec:ito} discusses the connection to It\^o formulas. Section~\ref{sec:analytical} records analytical consequences. Section~\ref{sec:mixed} treats mixed drivers. Section~\ref{sec:extensions} sketches extensions and open problems.

%=============================================================================
\section{Setup}
\label{sec:setup}
%=============================================================================

\paragraph{Probability space and filtrations.}
Throughout, $(\Omega, \mathcal{F}, P)$ denotes a complete probability space with a filtration $(\mathcal{F}_t)_{t \in [0,T]}$ satisfying the usual conditions (see, e.g., \cite{Rogers2000}): each $\mathcal{F}_t$ contains all $P$-null sets, and $\mathcal{F}_t = \mathcal{F}_t^+ := \bigcap_{s > t} \mathcal{F}_s$ (right-continuity). For a process~$X$, we write $(\mathcal{F}_t^X)$ for its natural filtration (augmented to satisfy the usual conditions). The predictable $\sigma$-algebra on $[0,T] \times \Omega$ is denoted~$\mathcal{P}$. The left-limit $\sigma$-algebra is $\mathcal{F}_{t-} := \sigma\!\bigl(\bigcup_{s < t} \mathcal{F}_s\bigr)$.

\paragraph{Function spaces.}
We write $L^2(\Omega)$ for square-integrable random variables with norm $\|F\|_{L^2} := \E[|F|^2]^{1/2}$. For a Hilbert space~$\mathcal{H}$, we denote its inner product by $\langle \cdot, \cdot \rangle_{\mathcal{H}}$ and norm by $\|\cdot\|_{\mathcal{H}}$. The domain of an operator~$A$ is $\dom(A)$; its adjoint (when it exists) is~$A^*$. We use $C_b^k(\R^d)$ for $k$-times continuously differentiable functions with bounded derivatives.

\paragraph{Operator conventions.}
The following notation, defined formally in Section~\ref{sec:OCD}, is used throughout: for a process~$X$ with energy space~$\mathcal{H}_X$, stochastic integral $\delta_X \colon \mathcal{H}_X \to L^2(\Omega)$, and operator--covariant derivative $D_X = \delta_X^*$, we write $\Pi_X$ for the predictable projection on~$\mathcal{H}_X$.

%=============================================================================
\section{The Operator--Covariant Derivative}
\label{sec:OCD}
%=============================================================================

This section introduces the main construction: a stochastic derivative defined intrinsically via adjointness to the stochastic integral.

\subsection{Energy spaces}

\begin{definition}[Energy Space]\label{def:energy}
Let $X = (X_t)_{t \in [0,T]}$ be a square-integrable process on $(\Omega, \mathcal{F}, P)$. An \emph{energy space} for~$X$ is a real separable Hilbert space~$\mathcal{H}_X$ of (equivalence classes of) processes $u = (u_t)_{t \in [0,T]}$, together with a continuous linear operator
\[
\delta_X \colon \mathcal{H}_X \to L^2(\Omega)
\]
called the \emph{stochastic integral} (or \emph{divergence}) with respect to~$X$.
\end{definition}

\begin{remark}[Bounded Integrals]\label{rem:bounded}
The continuity (boundedness) of $\delta_X$ is essential: it ensures via Riesz representation that the adjoint $D_X = \delta_X^*$ is a bounded operator on all of $L^2(\Omega)$. In much of the Malliavin calculus literature, the divergence/Skorokhod operator is treated as an unbounded, closed, densely defined operator~\cite{Nualart2006}. Our framework does not cover that setting. In the examples treated in this paper---It\^o integrals with respect to martingales and compensated Poisson processes---the stochastic integral is bounded (in fact isometric) on the energy space, so the bounded setting applies directly. Extending the framework to closed unbounded divergence operators, where $D_X$ would be defined only on a dense subspace of $L^2(\Omega)$, is an interesting direction for future work.
\end{remark}

\begin{example}[Martingale Energy Space]\label{ex:martingale_energy}
Let $X$ be a continuous square-integrable martingale with quadratic variation $\langle X \rangle_t$. The energy space is
\[
\mathcal{H}_X := \overline{\{\text{simple predictable processes}\}}^{\|\cdot\|_{\mathcal{H}_X}},
\]
where $\|u\|_{\mathcal{H}_X}^2 := \E\!\bigl[\int_0^T |u_t|^2\, d\langle X \rangle_t\bigr]$. The stochastic integral $\delta_X(u) := \int_0^T u_t\, dX_t$ is the It\^o integral, which extends by isometry to all of~$\mathcal{H}_X$.
\end{example}

\begin{definition}[Predictable Projection on Energy Space]\label{def:pred_proj}
Assume that $\mathcal{H}_X$ consists of jointly measurable processes and that the predictable processes form a closed subspace of~$\mathcal{H}_X$. The \emph{predictable projection} $\Pi_X \colon \mathcal{H}_X \to \mathcal{H}_X$ is the orthogonal projection onto this subspace.

The closedness assumption is not automatic for an arbitrary energy space; it must be verified. In all examples of this paper, the predictable subspace is either all of~$\mathcal{H}_X$ (for It\^o integrals with respect to martingales and compensated Poisson processes, where $\mathcal{H}_X$ consists entirely of predictable processes) or is closed by standard results on conditional expectation.

For the integral-form energy spaces used throughout this paper---specifically, those where $u_t \in L^2(\Omega)$ for each~$t$ and the $\mathcal{H}_X$ inner product takes the form $\langle u, v \rangle_{\mathcal{H}_X} = \E[\int_0^T u_t v_t\, d\mu_t]$ for some positive measure $\mu$ on $[0,T]$ (possibly random, as with $d\langle X\rangle_t$)---the predictable projection admits the pointwise representation
\[
(\Pi_X u)_t = \E[u_t \mid \mathcal{F}_{t-}] \quad \text{for } t \in (0,T], \qquad (\Pi_X u)_0 = \E[u_0 \mid \mathcal{F}_0],
\]
where $\mathcal{F}_{t-}$ is the left-limit $\sigma$-algebra defined in Section~\ref{sec:setup}. This is the standard predictable projection of optional processes; see, e.g., \cite[Chapter~VI]{Revuz1999}.
\end{definition}

\begin{remark}[When $\Pi_X$ Is Nontrivial]\label{rem:Pi_nontrivial}
For martingale integrals (Example~\ref{ex:martingale_energy}), the energy space~$\mathcal{H}_X$ consists of predictable processes, so $\Pi_X = \Id$. The projection $\Pi_X$ would become nontrivial in settings where $\mathcal{H}_X$ includes non-predictable elements. In the Skorokhod integral framework on Wiener space, for instance, the energy space is $L^2([0,T] \times \Omega)$ including anticipating processes; however, the Skorokhod integral is unbounded, so it falls outside the scope of the present paper (see Remark~\ref{rem:bounded}).
\end{remark}

\begin{remark}[Dependence on Filtration]\label{rem:filtration_dep}
The factorization~\eqref{eq:factorization} depends on three choices: the energy space~$\mathcal{H}_X$, the stochastic integral~$\delta_X$, and the filtration $(\mathcal{F}_t)$ (which determines the predictable projection~$\Pi_X$). The operator $D_X = \delta_X^*$ is determined by $\mathcal{H}_X$ and $\delta_X$ alone, but the factorization as a whole depends on the filtration through~$\Pi_X$. Different filtrations for the same process~$X$ will generally yield different factorizations. In all examples of this paper, we use the natural filtration of~$X$.
\end{remark}

\subsection{Definition by Riesz representation}

The following definition is the core of the paper.

\begin{definition}[Operator--Covariant Derivative]\label{def:OCD}
Let $X$ be a square-integrable process with energy space~$\mathcal{H}_X$ and stochastic integral $\delta_X \colon \mathcal{H}_X \to L^2(\Omega)$. For $F \in L^2(\Omega)$, the \emph{operator--covariant derivative} $D_X F \in \mathcal{H}_X$ is the unique element satisfying
\begin{equation}\label{eq:adjoint}
\E[F \cdot \delta_X(u)] \;=\; \langle D_X F,\, u \rangle_{\mathcal{H}_X}
\quad \text{for all } u \in \mathcal{H}_X.
\end{equation}
The operator $D_X \colon L^2(\Omega) \to \mathcal{H}_X$ is called the \emph{operator--covariant derivative}.
\end{definition}

\begin{remark}[Nature of the ``Derivative'']\label{rem:key_point}
The operator $D_X$ is not a derivative in any infinitesimal, pathwise, or variational sense, nor is it a derivation in the algebraic sense (i.e., it does not satisfy the Leibniz rule; see Section~\ref{ssec:leibniz}). It is the Riesz representer of the bounded linear functional $u \mapsto \E[F \cdot \delta_X(u)]$ on~$\mathcal{H}_X$. We retain the term ``derivative'' solely because of the adjoint relationship $D_X = \delta_X^*$, which mirrors the classical duality between differentiation and integration, and because $D_X$ recovers known stochastic derivatives in concrete settings. In the Brownian case, $D_B F$ recovers the conditional Malliavin derivative $\E[D_t^{\Mall} F \mid \mathcal{F}_t]$, which is a genuine derivative in the Sobolev sense; but this is a consequence, not part of the definition. See Section~\ref{ssec:leibniz} for the product structure of~$D_X$.
\end{remark}

\begin{remark}[Scope of the Domain]\label{rem:domain}
The adjoint $D_X = \delta_X^*$ is a bounded linear operator on all of $L^2(\Omega)$, by the Riesz representation theorem. However, the factorization~\eqref{eq:factorization} and representation~\eqref{eq:CO} hold only for $\mathcal{F}_T^X$-measurable~$F$, and require the additional hypotheses of Theorem~\ref{thm:CO} (centered integrals, isometry, representation property). For $F$ that is not representable, $D_X F$ is well defined as an element of~$\mathcal{H}_X$ and $\delta_X(\Pi_X D_X F)$ provides the best $L^2$-approximation of $F - \E[F]$ by stochastic integrals of~$X$; see Proposition~\ref{prop:GKW}.
\end{remark}

\begin{proposition}[Existence and Uniqueness]\label{prop:existence}
If $\delta_X \colon \mathcal{H}_X \to L^2(\Omega)$ is continuous, then for each $F \in L^2(\Omega)$, there exists a unique $D_X F \in \mathcal{H}_X$ satisfying~\eqref{eq:adjoint}.
\end{proposition}

\begin{proof}
The map $u \mapsto \E[F \cdot \delta_X(u)]$ is a bounded linear functional on~$\mathcal{H}_X$:
\[
|\E[F \cdot \delta_X(u)]| \;\le\; \|F\|_{L^2(\Omega)}\,\|\delta_X(u)\|_{L^2(\Omega)} \;\le\; C\,\|F\|_{L^2(\Omega)}\,\|u\|_{\mathcal{H}_X},
\]
where $C := \|\delta_X\|_{\op}$ is the operator norm. By the Riesz representation theorem, there exists a unique $D_X F \in \mathcal{H}_X$ with $\E[F \cdot \delta_X(u)] = \langle D_X F, u \rangle_{\mathcal{H}_X}$.
\end{proof}

\begin{proposition}[Basic Properties of $D_X$]\label{prop:properties}
The operator--covariant derivative $D_X \colon L^2(\Omega) \to \mathcal{H}_X$ satisfies:
\begin{enumerate}[label=\textup{(\roman*)}]
\item \textup{Linearity:} $D_X(\alpha F + \beta G) = \alpha\, D_X F + \beta\, D_X G$ for $\alpha, \beta \in \R$.
\item \textup{Boundedness:} $\|D_X F\|_{\mathcal{H}_X} \le C\,\|F\|_{L^2(\Omega)}$ where $C = \|\delta_X\|_{\op}$.
\item \textup{Adjointness:} $D_X = \delta_X^*$, i.e., $D_X$ is the Hilbert space adjoint of~$\delta_X$.
\end{enumerate}
\end{proposition}

\begin{proof}
Parts (i) and (iii) follow directly from~\eqref{eq:adjoint}. For~(ii), by Riesz representation:
\[
\|D_X F\|_{\mathcal{H}_X}
= \sup_{\|u\|_{\mathcal{H}_X} \le 1} |\langle D_X F, u \rangle_{\mathcal{H}_X}|
= \sup_{\|u\|_{\mathcal{H}_X} \le 1} |\E[F \cdot \delta_X(u)]|
\le C\,\|F\|_{L^2(\Omega)}. \qedhere
\]
\end{proof}

\subsection{Unified Clark--Ocone representation}

\begin{theorem}[Unified Clark--Ocone Representation]\label{thm:CO}
Let $X$ be a square-integrable process with natural filtration $(\mathcal{F}_t^X)_{t \in [0,T]}$, energy space~$\mathcal{H}_X$, stochastic integral~$\delta_X$, predictable projection~$\Pi_X$, and operator--covariant derivative $D_X = \delta_X^*$. Assume:
\begin{enumerate}[label=\textup{(\roman*)}]
\item \textup{Centered integrals:} $\E[\delta_X(u)] = 0$ for all $u \in \mathcal{H}_X$.
\item \textup{Isometry on predictable processes:} $\E[|\delta_X(u)|^2] = \|u\|_{\mathcal{H}_X}^2$ for predictable $u \in \mathcal{H}_X$.
\item \textup{Representation property:} Every $\mathcal{F}_T^X$-measurable $\tilde{F} \in L^2(\Omega)$ with $\E[\tilde{F}] = 0$ can be written as $\tilde{F} = \delta_X(v)$ for some predictable $v \in \mathcal{H}_X$.
\end{enumerate}
Then for all $\mathcal{F}_T^X$-measurable $F \in L^2(\Omega)$:
\begin{equation}\label{eq:CO_thm}
F \;=\; \E[F] + \delta_X(\Pi_X D_X F) \;=\; \E[F] + \int_0^T (\Pi_X D_X F)_t\, dX_t.
\end{equation}
Equivalently, the operator identity $(\Id - \E) = \delta_X \circ \Pi_X \circ D_X$ holds on $L^2(\Omega, \mathcal{F}_T^X, P)$.
\end{theorem}

\begin{proof}
By hypothesis~(iii), every centered $\tilde{F} := F - \E[F] \in L^2(\Omega)$ can be written as $\tilde{F} = \delta_X(v)$ for some predictable $v \in \mathcal{H}_X$. For any $u \in \mathcal{H}_X$, by adjointness~\eqref{eq:adjoint}:
\[
\langle D_X F, u \rangle_{\mathcal{H}_X}
= \E[F \cdot \delta_X(u)]
= \E[\tilde{F} \cdot \delta_X(u)] + \E[F] \cdot \E[\delta_X(u)]
= \E[\delta_X(v) \cdot \delta_X(u)],
\]
where the last equality uses $\tilde{F} = \delta_X(v)$ and hypothesis~(i). For predictable~$u$, hypothesis~(ii) (extended by polarization) gives $\E[\delta_X(v) \cdot \delta_X(u)] = \langle v, u \rangle_{\mathcal{H}_X}$. Thus
\[
\langle D_X F, u \rangle_{\mathcal{H}_X} = \langle v, u \rangle_{\mathcal{H}_X}
\quad \text{for all predictable } u.
\]
Since $\Pi_X$ is orthogonal projection onto predictable processes, $\langle \Pi_X D_X F, u \rangle_{\mathcal{H}_X} = \langle D_X F, u \rangle_{\mathcal{H}_X} = \langle v, u \rangle_{\mathcal{H}_X}$ for all predictable~$u$. As both $\Pi_X D_X F$ and $v$ are predictable, this implies $\Pi_X D_X F = v$, and thus $\tilde{F} = \delta_X(\Pi_X D_X F)$.
\end{proof}

\begin{remark}[Strength of the Representation Property]\label{rem:PRP}
Hypothesis~(iii) is the predictable representation property (PRP): it asserts that $\delta_X$ maps predictable integrands surjectively onto centered $\mathcal{F}_T^X$-measurable elements of~$L^2(\Omega)$. This is a strong structural condition that holds for Brownian motion and for compensated Poisson processes (each with respect to its own filtration), but fails for many natural models (e.g., fractional Brownian motion, mixed-filtration models). For a general L\'evy process with both continuous and jump components, the representation requires \emph{two} integrals (continuous martingale plus compensated Poisson random measure), which corresponds to the direct-sum framework of Section~\ref{sec:mixed}, not the single-integral PRP of hypothesis~(iii). The factorization $(\Id - \E) = \delta_X \circ \Pi_X \circ D_X$ requires PRP; without it, $D_X$ is a well-defined adjoint but the identity~\eqref{eq:CO_thm} need not hold. The representation integrand $\Pi_X D_X F$ is unique: the isometry (ii) implies that $\delta_X$ is injective on predictable processes, so $\delta_X(v_1) = \delta_X(v_2)$ with $v_1, v_2$ predictable implies $v_1 = v_2$ in~$\mathcal{H}_X$.
\end{remark}

\begin{remark}[When $\Pi_X = \Id$]\label{rem:Pi_id}
In all examples of this paper where Theorem~\ref{thm:CO} applies---Brownian motion, general continuous martingales, and compensated Poisson processes---the energy space~$\mathcal{H}_X$ consists entirely of predictable processes, so $\Pi_X = \Id$ and the factorization simplifies to $(\Id - \E) = \delta_X \circ D_X$. In these cases, $D_X F$ is already predictable and $\Pi_X$ plays no role. The general formulation retains~$\Pi_X$ for two reasons: first, it ensures structural correctness if the framework is applied in a setting where $\mathcal{H}_X$ admits non-predictable elements (as would occur in the unbounded Skorokhod extension noted in Section~\ref{sec:extensions}); second, even in the predictable case, the formula $\delta_X \circ \Pi_X \circ D_X$ makes explicit that the integrand is predictable, rather than relying on the reader to check that $D_X F$ happens to be so.
\end{remark}

\subsection{Derivations, cocycles, and correction terms}
\label{ssec:leibniz}

Classical differential calculus relies on the Leibniz rule $D(FG) = F \cdot DG + G \cdot DF$. In stochastic calculus, this rule typically fails: the It\^o formula includes a correction term $\tfrac{1}{2} f''(X)\,d\langle X \rangle$, and jump processes introduce additional corrections.

\begin{remark}[Leibniz Structure of $D_X$]\label{rem:leibniz}
The operator--covariant derivative $D_X = \delta_X^*$ is not a derivation in general.
\begin{enumerate}[label=\textup{(\roman*)}]
\item \emph{Zero quadratic variation:} If $X$ has zero quadratic variation and no jumps, then $D_X$ satisfies Leibniz. This is the classical deterministic case.
\item \emph{Continuous martingales:} For Brownian motion, the Malliavin derivative $D^{\Mall}$ satisfies a Leibniz rule, while the operator $D_B$ (which equals the predictable projection of $D^{\Mall}$; see Example~\ref{ex:brownian}) absorbs the conditional expectation. The cocycle or correction relating $D_X(FG)$ to $F \cdot D_X G + G \cdot D_X F$ is connected to the carr\'e du champ operator and the $\tfrac{1}{2}f''(X)\,d\langle X\rangle$ term of It\^o's formula; see~\cite{Nualart2006} for the precise relationship in the Malliavin setting.
\item \emph{Jump processes:} Additional correction terms arise from compensated sums over jumps, reflecting $f(X_t) - f(X_{t-}) - f'(X_{t-})\Delta X_t$, where $\Delta X_t := X_t - X_{t-}$ denotes the jump at time~$t$.
\end{enumerate}
\end{remark}

\begin{remark}[Why Factorization Survives]\label{rem:factorization_survives}
The factorization $(\Id - \E) = \delta_X \circ \Pi_X \circ D_X$ does not require $D_X$ to satisfy Leibniz. The identity is a statement about representations, not about product rules:
\begin{enumerate}[label=\textup{(\roman*)}]
\item The Clark--Ocone representation~\eqref{eq:CO_thm} holds under the hypotheses of Theorem~\ref{thm:CO}, regardless of whether $D_X$ satisfies any product rule.
\item The It\^o formula, which involves products and hence correction terms, arises separately by differentiating conditional expectations along paths of~$X$ (see Section~\ref{sec:ito}). The correction terms appear at that stage, not in the factorization itself.
\end{enumerate}
\end{remark}

%=============================================================================
\section{Concrete Examples}
\label{sec:examples}
%=============================================================================

We verify that the operator--covariant framework recovers classical derivatives as special cases and provide a substantive non-Gaussian example.

\subsection{Brownian motion}

\begin{example}[Brownian Case]\label{ex:brownian}
Let $X = B$ be standard Brownian motion with its natural filtration. Then:
\begin{enumerate}[label=\textup{(\roman*)}]
\item The energy space~$\mathcal{H}_B$ is the space of predictable processes $u = (u_t)_{t \in [0,T]}$ with $\|u\|_{\mathcal{H}_B}^2 := \E\!\bigl[\int_0^T |u_t|^2\, dt\bigr] < \infty$.
\item The stochastic integral $\delta_B(u) = \int_0^T u_t\, dB_t$ is the It\^o integral, satisfying the It\^o isometry $\E[|\delta_B(u)|^2] = \|u\|_{\mathcal{H}_B}^2$.
\item Since $\mathcal{H}_B$ consists of predictable processes, $\Pi_B = \Id$, and by the martingale representation theorem, the representation~\eqref{eq:CO_thm} becomes:
\[
F = \E[F] + \int_0^T (D_B F)_t\, dB_t.
\]
\end{enumerate}
For $F \in \mathbb{D}^{1,2}$ (the Malliavin--Sobolev space of differentiable random variables; see~\cite{Nualart2006}), the classical Clark--Ocone formula gives $F - \E[F] = \int_0^T \E[D_t^{\Mall} F \mid \mathcal{F}_t]\, dB_t$, where $D^{\Mall}$ is the Malliavin derivative. Comparing, we see $(D_B F)_t = \E[D_t^{\Mall} F \mid \mathcal{F}_t]$ for such~$F$. The operator $D_B$ is defined on all of $L^2(\Omega)$ via adjointness (though representationally meaningful only for $\mathcal{F}_T^B$-measurable~$F$), whereas $D^{\Mall}$ is defined only on the dense subspace $\mathbb{D}^{1,2}$.
\end{example}

\begin{remark}[Comparison with Malliavin Theory]\label{rem:malliavin}
In classical Malliavin calculus, one defines a derivative $D^{\Mall}$ on a dense subspace $\mathbb{D}^{1,2} \subset L^2(\Omega)$ of differentiable functionals. The Clark--Ocone formula then states $F - \E[F] = \int_0^T \E[D_t^{\Mall} F \mid \mathcal{F}_t]\, dB_t$ for $F \in \mathbb{D}^{1,2}$. Our approach reverses this: $D_B$ is defined on all of $L^2(\Omega)$ via adjointness, without requiring differentiability. For $F \in \mathbb{D}^{1,2}$, we have $(D_B F)_t = \E[D_t^{\Mall} F \mid \mathcal{F}_t]$---that is, $D_B F$ is the predictable projection of the Malliavin derivative.

The trade-off is clear: $D^{\Mall}$ carries richer pointwise information (it is a random field indexed by~$t$ before conditioning), while $D_B$ directly provides the predictable integrand needed for representation. Neither subsumes the other; they answer different questions.
\end{remark}

\subsection{General continuous martingales}

\begin{example}[General Continuous Martingale]\label{ex:general_mart}
Let $X$ be any continuous square-integrable martingale with quadratic variation $\langle X \rangle$ and the martingale representation property. The energy space~$\mathcal{H}_X$ consists of predictable processes $u = (u_t)_{t \in [0,T]}$ with $\E\!\bigl[\int_0^T u_t^2\, d\langle X \rangle_t\bigr] < \infty$. Since $\mathcal{H}_X$ already consists of predictable processes, $\Pi_X = \Id$ and the operator--covariant derivative $D_X$ satisfies
\[
F = \E[F] + \int_0^T (D_X F)_t\, dX_t
\]
for all $\mathcal{F}_T^X$-measurable $F \in L^2(\Omega)$.
\end{example}

\subsection{Compensated Poisson process}
\label{ssec:poisson}

This is a genuinely non-Gaussian example where we can compute the operator--covariant derivative explicitly.

Let $N = (N_t)_{t \in [0,T]}$ be a Poisson process with intensity $\lambda > 0$, and let $\tilde{N}_t := N_t - \lambda t$ be the compensated Poisson process. The natural filtration is $(\mathcal{F}_t^N)$.

\begin{example}[Compensated Poisson Process]\label{ex:poisson}
The energy space is
\[
\mathcal{H}_{\tilde{N}} := \bigl\{ u = (u_t)_{t \in [0,T]} : u \text{ is predictable},\; \|u\|_{\mathcal{H}_{\tilde{N}}}^2 := \lambda\, \E\!\bigl[\textstyle\int_0^T |u_t|^2\, dt\bigr] < \infty \bigr\},
\]
with stochastic integral $\delta_{\tilde{N}}(u) := \int_0^T u_t\, d\tilde{N}_t$. The isometry
\[
\E\!\bigl[|\delta_{\tilde{N}}(u)|^2\bigr] = \lambda\, \E\!\bigl[\textstyle\int_0^T |u_t|^2\, dt\bigr] = \|u\|_{\mathcal{H}_{\tilde{N}}}^2
\]
holds for predictable~$u$.
\end{example}

\begin{proposition}[Representation Property for $\tilde{N}$]\label{prop:poisson_rep}
The compensated Poisson process $\tilde{N}$ has the representation property: every $\mathcal{F}_T^N$-measurable $F \in L^2(\Omega)$ with $\E[F] = 0$ can be written as $F = \int_0^T v_t\, d\tilde{N}_t$ for some predictable $v \in \mathcal{H}_{\tilde{N}}$.
\end{proposition}

\begin{proof}
This is a standard result; see, e.g., \cite[Chapter~12]{Kallenberg2002} or \cite[Section~1.5]{DiNunno2009}. The compensated Poisson process is a martingale with the predictable representation property.
\end{proof}

Since $\mathcal{H}_{\tilde{N}}$ consists of predictable processes, $\Pi_{\tilde{N}} = \Id$. The hypotheses of Theorem~\ref{thm:CO} are satisfied: (i)~$\E[\delta_{\tilde{N}}(u)] = 0$ because the compensated Poisson integral of a predictable process is a martingale; (ii)~the isometry holds by Example~\ref{ex:poisson}; (iii)~the representation property holds by Proposition~\ref{prop:poisson_rep}. Theorem~\ref{thm:CO} therefore gives:

\begin{corollary}\label{cor:poisson_CO}
For all $\mathcal{F}_T^N$-measurable $F \in L^2(\Omega)$:
\[
F = \E[F] + \int_0^T (D_{\tilde{N}} F)_t\, d\tilde{N}_t.
\]
\end{corollary}

We now identify $D_{\tilde{N}}$ explicitly. For Poisson functionals, the natural difference operator is the ``add-a-point'' operator of the discrete Malliavin calculus~\cite{Privault2009}. Identifying each sample path $\omega$ of the Poisson process with its set of jump times $\{T_1, T_2, \ldots\} \subset [0,T]$, we write $\omega + \varepsilon_t$ for the configuration obtained by adding one jump at time~$t$. The add-a-point operator is then
\begin{equation}\label{eq:diff_op}
(\mathcal{D}_t^+ G)(\omega) := G(\omega + \varepsilon_t) - G(\omega).
\end{equation}

\begin{proposition}[Explicit Form of $D_{\tilde{N}}$]\label{prop:poisson_explicit}
For $\mathcal{F}_T^N$-measurable $F \in L^2(\Omega)$, the operator--covariant derivative is
\[
(D_{\tilde{N}} F)_t = \E[\mathcal{D}_t^+ F \mid \mathcal{F}_{t-}]
\]
$d t \otimes dP$-a.e., where $\mathcal{D}_t^+$ is the add-a-point operator~\eqref{eq:diff_op}. In particular, $D_{\tilde{N}} F$ is the predictable projection of the discrete gradient. This identification is well known in the discrete Malliavin calculus~\cite{Privault2009,DiNunno2009}; what is new here is not the operator itself, but its characterization as an instance of the general adjoint construction $D_X = \delta_X^*$.
\end{proposition}

\begin{proof}
We verify the adjoint relation~\eqref{eq:adjoint}. For predictable $u \in \mathcal{H}_{\tilde{N}}$, the integration-by-parts formula for the Poisson process (see \cite[Proposition~1.9]{Privault2009} or \cite[Theorem~12.5]{DiNunno2009}) gives:
\[
\E\!\Bigl[F \cdot \int_0^T u_t\, d\tilde{N}_t\Bigr]
= \lambda\, \E\!\Bigl[\int_0^T u_t\, \mathcal{D}_t^+ F\, dt\Bigr].
\]
Since $u_t$ is predictable (hence $\mathcal{F}_{t-}$-measurable), the tower property yields
\[
\E[u_t\, \mathcal{D}_t^+ F] = \E\bigl[u_t\, \E[\mathcal{D}_t^+ F \mid \mathcal{F}_{t-}]\bigr]
\]
for each~$t$. Thus $\E[F \cdot \delta_{\tilde{N}}(u)] = \lambda\, \E\!\bigl[\int_0^T u_t\, \E[\mathcal{D}_t^+ F \mid \mathcal{F}_{t-}]\, dt\bigr]$. Since $\langle g, u \rangle_{\mathcal{H}_{\tilde{N}}} = \lambda\,\E\!\bigl[\int_0^T g_t\, u_t\, dt\bigr]$ for predictable~$g$ and~$u$, this reads
\[
\E[F \cdot \delta_{\tilde{N}}(u)] = \langle \E[\mathcal{D}_\cdot^+ F \mid \mathcal{F}_{\cdot-}],\, u \rangle_{\mathcal{H}_{\tilde{N}}}.
\]
By uniqueness in Definition~\ref{def:OCD}, $(D_{\tilde{N}} F)_t = \E[\mathcal{D}_t^+ F \mid \mathcal{F}_{t-}]$.
\end{proof}

\begin{example}[Explicit Computation]\label{ex:poisson_compute}
Let $F = N_T^2$. Then $\E[F] = \lambda T + \lambda^2 T^2$. The add-a-point operator gives
\[
\mathcal{D}_t^+ (N_T^2) = (N_T + 1)^2 - N_T^2 = 2N_T + 1.
\]
Since $N_T - N_{t-}$ is independent of $\mathcal{F}_{t-}$ (by the independent increments property, noting that $\Delta N_t = 0$ a.s.\ for any fixed~$t$):
\[
\E[\mathcal{D}_t^+ (N_T^2) \mid \mathcal{F}_{t-}]
= 2\E[N_T \mid \mathcal{F}_{t-}] + 1
= 2\bigl(N_{t-} + \lambda(T - t)\bigr) + 1
= 2N_{t-} + 2\lambda(T - t) + 1.
\]
Thus $(D_{\tilde{N}} F)_t = 2N_{t-} + 2\lambda(T-t) + 1$, and the representation reads
\[
N_T^2 = \E[N_T^2] + \int_0^T \bigl(2N_{t-} + 2\lambda(T-t) + 1\bigr)\, d\tilde{N}_t.
\]
One can verify this directly by expanding the stochastic integral and using properties of the Poisson process.
\end{example}

\subsection{Scope boundary: Volterra processes}
\label{ssec:volterra}

We briefly discuss Volterra processes to delineate where the framework does \emph{not} apply in its current form.

\begin{definition}[Volterra Process]\label{def:volterra}
A Volterra Gaussian process is $X_t = \int_0^t K(t,s)\, dW_s$, where $W$ is Brownian motion and $K \colon \{(t,s) : 0 \le s \le t \le T\} \to \R$ is a square-integrable kernel.
\end{definition}

The most important special case is fractional Brownian motion with Hurst parameter $H \in (0,1)$, which is the centered Gaussian process with covariance $R_H(s,t) = \tfrac{1}{2}(s^{2H} + t^{2H} - |t-s|^{2H})$.

\begin{remark}[Volterra Case: Approximation without Representation]\label{rem:volterra}
For a Volterra process $X_t = \int_0^t K(t,s)\, dW_s$, one can define an energy space and stochastic integral (see, e.g., \cite{Nualart2006}), and the adjoint $D_X := \delta_X^*$ is well defined whenever $\delta_X$ is continuous. However, Volterra processes are generally not martingales (since $K(t,s)$ depends on~$t$), so the hypotheses of Theorem~\ref{thm:CO}---particularly the representation property---are not satisfied, and the exact factorization $(\Id - \E) = \delta_X \circ \Pi_X \circ D_X$ does not hold.

This is a genuine limitation: the exact representation requires PRP, which excludes processes such as fractional Brownian motion with $H \neq 1/2$. Clark--Ocone representations for Volterra processes in the literature are derived using Malliavin calculus on the underlying Brownian motion~$W$, with the kernel~$K$ mediating the relationship~\cite{Nualart2006}.

Nevertheless, if the energy space and stochastic integral for~$X$ satisfy hypotheses~(i)--(ii) of Theorem~\ref{thm:CO} (centered integrals and isometry on predictable processes), then by Proposition~\ref{prop:GKW} the term $\delta_X(\Pi_X D_X F)$ provides the optimal $L^2$-approximation of $F - \E[F]$ by stochastic integrals of~$X$. Thus, even when the exact representation fails, the operator--covariant derivative identifies the best mean-square approximation. Extending the framework to obtain exact representation theorems for non-martingale Volterra processes remains an open problem (Section~\ref{sec:extensions}).
\end{remark}

%=============================================================================
\section{Connection to It\^o Formulas}
\label{sec:ito}
%=============================================================================

The Clark--Ocone representation expresses a random variable as an integral. A natural question is: how does this relate to the classical It\^o formula? The connection runs through conditional expectations and backward equations. We give a brief sketch; a full development is beyond the scope of this paper.

For a diffusion $dX_t = b(X_t)\,dt + \sigma(X_t)\,dB_t$ with $X_0 = x_0$, where $b, \sigma$ satisfy standard Lipschitz and linear growth conditions ensuring a unique strong solution (so that $\mathcal{F}_t^X = \mathcal{F}_t^B$ up to null sets), and $F = f(X_T)$ where $f \in C_b^2(\R^d)$, define $u(t,x) := \E[f(X_T) \mid X_t = x]$. The function~$u$ satisfies the backward Kolmogorov equation $\partial_t u + \mathcal{L}u = 0$ on $[0,T) \times \R^d$ with $u(T,x) = f(x)$, where $\mathcal{L} := \tfrac{1}{2}\sum_{i,j}(\sigma\sigma^\top)_{ij}(x)\,\partial_{x_i}\partial_{x_j} + \sum_i b_i(x)\,\partial_{x_i}$ is the infinitesimal generator of~$X$ (see, e.g., \cite{Oksendal2003}). Since $M_t := u(t, X_t)$ is a martingale, applying It\^o's formula to~$M_t$ and comparing with the Clark--Ocone representation $F = \E[F] + \int_0^T (D_B F)_t\, dB_t$ identifies the integrand $(D_B F)_t$ with $\nabla u(t, X_t)^\top \sigma(X_t)$.

The same principle extends to path-dependent functionals via the Dupire derivatives of functional It\^o calculus~\cite{Cont2013,Dupire2009}: for $F = U(T, X_{[0,T]})$ with sufficiently regular non-anticipative functional~$U$, the operator~$D_X$ provides an intrinsic characterization of the integrand via adjointness, while the vertical derivative~$\nabla_\omega U$ characterizes it via path-level differentiation.

The factorization framework does not produce new It\^o formulas; rather, it shows that the integrand in the representation formula can be computed either by the adjoint route ($D_X$) or by classical differentiation along paths, and these must agree. More broadly, the factorization~\eqref{eq:factorization} is purely representation-theoretic: it provides a natural decomposition of $L^2$ random variables into stochastic integrals (relative to the chosen energy space and filtration), but does not by itself yield pathwise calculus, infinitesimal expansions, or chain rules. These differential identities arise only after imposing additional analytical structure (PDE regularity, path differentiability, etc.).

%=============================================================================
\section{Analytical Consequences}
\label{sec:analytical}
%=============================================================================

\begin{proposition}[Variance/Energy Identity]\label{prop:variance}
Under the hypotheses of Theorem~\ref{thm:CO}, for $\mathcal{F}_T^X$-measurable $F \in L^2(\Omega)$:
\[
\Var(F) = \|\Pi_X D_X F\|_{\mathcal{H}_X}^2.
\]
For continuous martingales with quadratic variation, this becomes
\[
\Var(F) = \E\!\bigl[\textstyle\int_0^T |(\Pi_X D_X F)_t|^2\, d\langle X \rangle_t\bigr].
\]
\end{proposition}

\begin{proof}
Apply the isometry to~\eqref{eq:CO_thm}:
\[
\Var(F) = \E[(F - \E[F])^2] = \E[|\delta_X(\Pi_X D_X F)|^2] = \|\Pi_X D_X F\|_{\mathcal{H}_X}^2. \qedhere
\]
\end{proof}

\begin{remark}[Energy Interpretation]\label{rem:energy}
Under the hypotheses of Theorem~\ref{thm:CO}, the integrand $\Pi_X D_X F$ is the unique predictable process representing $F - \E[F]$ as a stochastic integral, and its energy norm equals the variance of~$F$. This provides an operator-theoretic interpretation of the classical result that the variance of a martingale representation equals the $L^2$-norm of the integrand.
\end{remark}

\subsection{Beyond PRP: the optimal approximation property}

The definition of $D_X$ does not require PRP; only the factorization identity does. A natural question is: when PRP fails, what does $\delta_X(\Pi_X D_X F)$ represent? The following proposition shows it is the best $L^2$-approximation of $F - \E[F]$ by stochastic integrals of predictable processes---the Galtchouk--Kunita--Watanabe (GKW) projection.

\begin{proposition}[GKW Projection]\label{prop:GKW}
Let $X$ be a square-integrable process with energy space~$\mathcal{H}_X$ and bounded stochastic integral~$\delta_X$ satisfying hypotheses~(i)--(ii) of Theorem~\ref{thm:CO} (centered integrals and isometry on predictable processes). Let $\mathcal{S}_X := \overline{\{\delta_X(u) : u \in \mathcal{H}_X \text{ predictable}\}}^{L^2}$ be the stable subspace generated by predictable integrals of~$X$. Then for any $F \in L^2(\Omega)$:
\[
\delta_X(\Pi_X D_X F) = \mathrm{proj}_{\mathcal{S}_X}(F - \E[F]),
\]
where $\mathrm{proj}_{\mathcal{S}_X}$ denotes orthogonal projection onto~$\mathcal{S}_X$ in $L^2(\Omega)$.
\end{proposition}

\begin{proof}
The element $\delta_X(\Pi_X D_X F) \in \mathcal{S}_X$ since $\Pi_X D_X F$ is predictable. It remains to show orthogonality: for any predictable $u \in \mathcal{H}_X$,
\begin{align*}
\E\bigl[(F - \E[F] - \delta_X(\Pi_X D_X F)) \cdot \delta_X(u)\bigr]
&= \E[F \cdot \delta_X(u)] - \E[\delta_X(\Pi_X D_X F) \cdot \delta_X(u)] \\
&= \langle D_X F, u \rangle_{\mathcal{H}_X} - \langle \Pi_X D_X F, u \rangle_{\mathcal{H}_X} \\
&= \langle D_X F - \Pi_X D_X F, u \rangle_{\mathcal{H}_X} = 0,
\end{align*}
where the second equality uses adjointness~\eqref{eq:adjoint} and the isometry (hypothesis~(ii) extended by polarization), and the last uses the fact that $D_X F - \Pi_X D_X F$ is orthogonal to all predictable elements.
\end{proof}

\begin{remark}[Interpretation]\label{rem:GKW_interp}
Proposition~\ref{prop:GKW} shows that $D_X$ is meaningful even when PRP fails. When the representation property holds, $\mathcal{S}_X = \{G \in L^2(\mathcal{F}_T^X) : \E[G] = 0\}$ and the GKW projection recovers the full representation $F - \E[F] = \delta_X(\Pi_X D_X F)$. When PRP fails, $\delta_X(\Pi_X D_X F)$ is the variance-minimizing approximation of $F - \E[F]$ by stochastic integrals of~$X$. In mathematical finance, this is precisely the optimal hedging strategy for the claim~$F$ using the tradeable asset~$X$ under mean-square hedging~\cite{Revuz1999}.
\end{remark}

%=============================================================================
\section{Mixed Drivers}
\label{sec:mixed}
%=============================================================================

Many applications involve processes driven by multiple noise sources---for instance, a diffusion plus jumps, or Brownian motion plus an independent Poisson process. When the noise components are independent, the operator--covariant framework handles such settings naturally via direct sums of energy spaces. The independence assumption is essential: it ensures orthogonality of the component energy spaces and additivity of the adjoint decomposition.

\subsection{Brownian--Poisson mixture}
\label{ssec:BP}

Let $B$ be standard Brownian motion and $\tilde{N}$ an independent compensated Poisson process with intensity~$\lambda$ (as in Section~\ref{ssec:poisson}). For $X_t := \alpha B_t + \beta \tilde{N}_t$ with $\alpha, \beta > 0$:

\begin{definition}[Mixed Energy Space]\label{def:mixed_energy}
$\mathcal{H}_X := \mathcal{H}_B \oplus \mathcal{H}_{\tilde{N}}$ with inner product $\langle (u, v), (u', v') \rangle_{\mathcal{H}_X} := \langle u, u' \rangle_{\mathcal{H}_B} + \langle v, v' \rangle_{\mathcal{H}_{\tilde{N}}}$.
\end{definition}

\begin{definition}[Mixed Divergence]\label{def:mixed_div}
For $(u, v) \in \mathcal{H}_X$, define $\delta_X(u,v) := \alpha\, \delta_B(u) + \beta\, \delta_{\tilde{N}}(v)$.
\end{definition}

The operator--covariant derivative $D_X := \delta_X^*$ decomposes as $D_X F = (\alpha\, D_B F,\, \beta\, D_{\tilde{N}} F)$, as follows from the adjoint relation: for $(u,v) \in \mathcal{H}_X$, $\E[F \cdot \delta_X(u,v)] = \alpha\,\E[F \cdot \delta_B(u)] + \beta\,\E[F \cdot \delta_{\tilde{N}}(v)] = \langle (\alpha\, D_B F,\, \beta\, D_{\tilde{N}} F), (u,v) \rangle_{\mathcal{H}_X}$.

Since $B$ and $\tilde{N}$ are independent and both have the representation property with respect to their own filtrations, the pair $(B, \tilde{N})$ has the representation property with respect to $\mathcal{F}_T^B \vee \mathcal{F}_T^N$ (see, e.g., \cite[Theorem~18.12]{Kallenberg2002}). If $F$ is measurable with respect to this joint filtration, the factorization~\eqref{eq:factorization} yields:
\[
F = \E[F] + \alpha \int_0^T (D_B F)_t\, dB_t + \beta \int_0^T (D_{\tilde{N}} F)_t\, d\tilde{N}_t.
\]

\subsection{General square-integrable L\'evy drivers}

We briefly discuss a more general non-Gaussian setting to illustrate the scope of the framework. Let $L$ be a centered pure-jump L\'evy process on $\R$ (i.e., $\E[L_t] = 0$) with L\'evy measure $\nu$ satisfying $\int_{\R} |z|^2\, \nu(dz) < \infty$ (finite second moment). The process $L$ has the representation $L_t = \int_0^t \int_{\R} z\, \widetilde{\mathcal{N}}(ds, dz)$, where $\mathcal{N}((0,t] \times A) := \#\{s \in (0,t] : \Delta L_s \in A\}$ is the Poisson random measure, $\widetilde{\mathcal{N}}(dt, dz) := \mathcal{N}(dt,dz) - dt\,\nu(dz)$ is the compensated measure, and $\Delta L_s := L_s - L_{s-}$ denotes the jump at time~$s$.

The energy space for~$L$ is $\mathcal{H}_L := L^2([0,T] \times \R,\, dt \otimes \nu)$ consisting of predictable integrands $v = v(t,z)$, with stochastic integral $\delta_L(v) := \int_0^T \int_{\R} v(t,z)\, \widetilde{\mathcal{N}}(dt, dz)$. The It\^o isometry
\[
\E[|\delta_L(v)|^2] = \E\!\Bigl[\int_0^T \int_{\R} |v(t,z)|^2\, \nu(dz)\, dt\Bigr] = \|v\|_{\mathcal{H}_L}^2
\]
holds for predictable~$v$, so $\delta_L$ is bounded and the operator--covariant derivative $D_L := \delta_L^*$ is well defined on $L^2(\Omega)$ by Remark~\ref{rem:bounded}. For a mixed process $X_t = \alpha B_t + \beta L_t$ with $B$ an independent standard Brownian motion and $\alpha, \beta > 0$, one takes $\mathcal{H}_X := \mathcal{H}_B \oplus \mathcal{H}_L$ with $\delta_X(u,v) := \alpha\, \delta_B(u) + \beta\, \delta_L(v)$, and $D_X := \delta_X^*$ decomposes as $D_X F = (\alpha\, D_B F,\, \beta\, D_L F)$ by the same argument as in Section~\ref{ssec:BP}.

\begin{remark}[Infinite-variance processes]\label{rem:infinite_var}
Symmetric $\gamma$-stable L\'evy processes with $\gamma \in (0,2)$ have infinite second moments ($\int |z|^2\, \nu(dz) = \infty$), so they fall outside the $L^2$ framework of this paper. Extending the factorization to such processes---for instance via $L^p$ or Besov-type energy spaces---is an interesting open problem; see Section~\ref{sec:extensions}.
\end{remark}

%=============================================================================
\section{Extensions and Open Problems}
\label{sec:extensions}
%=============================================================================

We note several directions for future work:

\begin{enumerate}[label=(\arabic*),leftmargin=2em]
\item \emph{Unbounded divergence operators.} The present framework requires $\delta_X$ to be a bounded linear map into $L^2(\Omega)$ (Remark~\ref{rem:bounded}). In many settings---notably the Skorokhod integral on Wiener space---the divergence is an unbounded, closed, densely defined operator. Extending the factorization to such settings, where $D_X$ would be defined only on a dense subspace of $L^2(\Omega)$, is a natural next step. Related questions include whether the resulting operator is closable in a natural topology, and whether it generates a Dirichlet form connecting the factorization framework to the energy measure theory of~\cite{Nualart2006}.

\item \emph{Infinite-variance and stable processes.} Symmetric $\gamma$-stable L\'evy processes ($\gamma < 2$) have infinite second moments and fall outside the $L^2$ theory (Remark~\ref{rem:infinite_var}). Extending the framework via $L^p$ or Besov-type energy spaces would bring important non-Gaussian processes into scope.

\item \emph{Volterra and non-semimartingale processes.} The main limitation of Theorem~\ref{thm:CO} is the representation property, which fails for non-martingale Volterra processes like fractional Brownian motion with $H \neq 1/2$. Extending the factorization to such processes---perhaps by working with the underlying Brownian motion and transferring the representation via the Volterra kernel---is an important open problem.

\item \emph{Fractional regimes $H < 1/2$.} Domain characterizations of $D_X$ when the energy space has a more complex structure (e.g., fractional Sobolev spaces for rough processes).

\item \emph{Jump processes: further examples.} Systematic development of the factorization for L\'evy processes beyond the Poisson case, including infinite-activity processes with finite second moments.

\item \emph{Mild SPDEs.} Adjoint derivatives for mild solutions of stochastic PDEs in Banach spaces.

\item \emph{Higher-order objects.} Second-order operators from iterated applications of $D_X$ may connect to curvature-type quantities; however, no geometric structure (connection, metric compatibility, manifold) is defined in the present framework, so such connections remain heuristic and speculative.
\end{enumerate}

%=============================================================================
\section*{Acknowledgments}
%=============================================================================

AI tools were used as interactive assistants for drafting and mathematical exploration. The author is responsible for all content.

%=============================================================================

\end{document}